\documentclass{article}
\usepackage{graphicx}
\usepackage{amsmath}
\usepackage{amsfonts}
\usepackage{amssymb}
\newtheorem{theorem}{Theorem}

\newtheorem{corollary}[theorem]{Corollary}

\newtheorem{definition}[theorem]{Definition}

\newtheorem{proposition}[theorem]{Proposition}

\newenvironment{proof}[1][Proof]{\textbf{#1.} }{\ \rule{0.5em}{0.5em}}

\begin{document}

\title{Morse theory, Milnor fibers and minimality of hyperplane arrangements}
\author{Richard Randell\\The University of Iowa}
\maketitle
\begin{abstract}
Through the study of Morse theory on the associated Milnor fiber, we show that
complex hyperplane arrangement complements are minimal. \ That is, the
complement of any complex hyperplane arrangement has the homotopy type of a CW
complex in which the number of p-cells equals the p-th betti number.
\ Combining this result with recent work of Papadima and Suciu, one obtains a
characterization of when arrangement complements are Eilenberg-MacLane spaces.
\end{abstract}

\section{Introduction}

\footnotetext{1991 Mathematics Subject Classification. \ Primary 52C35, 55Q52;
Secondary 14M12, 32S22
\par
Key words and phrases. \ hyperplane arrangement, milnor fiber, morse
theory}Let $\mathcal{A}$ be a central arrangement of hyperplanes in
$\mathbb{C}^{\ell}$, by which is meant a finite collection $\left\{
H_{1},\ldots,H_{n}\right\}  $ where $H_{i}=\alpha_{i}^{-1}(0)$ and each
$\alpha_{i}$ is a linear homogeneous form in the variables $(z_{1}%
,\ldots,z_{\ell})$. \ We call $\mathcal{A}$ an $\ell$\emph{-arrangement}. \ We
let $M$ be the complement of the union of the hyperplanes
\[
M=\mathbb{C}^{\ell}\setminus\cup H_{i}%
\]
and we consider the Milnor fibration $Q:M\rightarrow\mathbb{C}^{\ast
}=\mathbb{C\setminus}0$ , with $Q=\prod\alpha_{i}.$ This map is actually a
fiber bundle with fiber $F=Q^{-1}(1)$ \cite{Milnor1}, called the Milnor fiber.
\ In this paper we will study $F$ via Morse theory, following work of L\^{e}
\cite{Le1}. \ The Milnor fiber also appears as an n-fold cover of the
complement of the projectivized $M^{\ast}=M/\mathbb{C}^{\ast}$, where the
action is the restriction of the Hopf action. \ The Hopf action is free and
one sees that the associated bundle is trivial over $M^{\ast}$ so that $M\cong
M^{\ast}\times\mathbb{C}^{\ast}$. \ The Morse functions we consider on $F$ are
actually equivariant with respect to the covering transformations of the
n-fold cover (these covering transformations are generated by the monodromy of
the Milnor fibration), thus inducing Morse functions on $M^{\ast}$, hence a CW
complex structure on the homotopy type of $M^{\ast}$. \ Since $M$ has the
homotopy type of $M^{\ast}\times S^{1}$ one gets a CW-complex structure on the
homotopy type of $M$ by using the product structure. \ These CW structures on
$M$ and $M^{\ast}$ have the smallest possible number of cells in each
dimension, showing that central and affine arrangement complements are minimal
in the sense of Papadima-Suciu \cite{PapaSucmin}. \ From their work a
criterion for when the arrangement complement is an Eilenberg-MacLane
$K(\pi,1)$-space follows immediately (Corollary \ref{iff}).

The number of cells for all these CW complexes is determined combinatorially
by the intersection lattice of the arrangement (see \cite{OrlikTeraoBook} for
general background on hyperplane arrangements). \ However, for arrangements in
different path components of the configuration space, one may have varying
topology. \ In our context this would be reflected via the attaching maps in
the CW complexes. \ In a paper in preparation we study these attaching maps,
as well as the behavior of these CW complexes in lattice isotopic families and
in degenerating families.

Papers of Orlik and Terao \cite{OTMorse} and Dimca \cite{DimcaMilnor} consider
Morse theory and the Milnor fiber. \ In particular, Dimca also shows that
arrangement complements are minimal. \ Our construction follows L\^{e}
\cite{Le1} in using linear sections and readily yields minimality and explicit computation.

I am pleased to thank Alexandru Dimca, Michael Falk and Alexander Suciu for
helpful conversations and comments. \ In particular, A. Dimca pointed out that
the equivariant Morse function on $F$ induces a Morse function on $M^{\ast}$,
thereby giving a more straightforward proof of Theorem 3.

\section{Morse theory on the Milnor fiber}

We follow L\^{e} \cite{Le1} throughout this section. \ By a CW structure on a
space $X$, we shall mean a CW complex homotopy equivalent to $X$. \ In his
work, L\^{e} treats a general hypersurface singularity. \ For a singularity
given by the defining function $Q$ of a hyperplane arrangement the situation
is especially pleasant. \ The basic idea is to obtain a CW structure on the
Milnor fiber by taking successive generic linear sections and noting that one
goes from a particular section to a section of (complex) dimension one greater
by attaching cells of known dimension according to a Morse function which is
simply the distance from a suitably chosen hyperplane section of the Milnor
fiber. \ The result from L\^{e} is

\begin{theorem}
(\cite{Le1}) For almost all linear forms $L$, the real-valued function
$\left|  L\right|  _{F}$ (which measures the distance from a point to the
hyperplane $L=0$) is a Morse function (relative to the subset $F\cap\{L=0\})
$. \ At each critical point this function has index $\ell-1$.
\end{theorem}

By a linear change of coordinates we may assume $L=z_{\ell}$, and we will
henceforth do so. \ As in \cite{Le1}, we define the polar curve $\Gamma
=\{\frac{\partial Q}{\partial z_{1}}=\frac{\partial Q}{\partial z_{2}%
}=...=\frac{\partial Q}{\partial z_{\ell-1}}=0\}\cap\{Q\neq0\}$.

\begin{proposition}
In the setting of this paper, $\left|  z_{\ell}\right|  _{F}$ is a relative
Morse function provided that the linear subspace $\left\{  z_{\ell}=0\right\}
\subset\mathbb{C}^{\ell}$ contains no intersection of the hyperplanes of the
arrangement other than the origin, and that the set $\Gamma$ is a reduced
(possibly empty) curve.
\end{proposition}

\begin{proof}
We check directly that L\^{e}'s genericity conditions are satisfied for such
$L=z_{\ell}$. \ First one stratifies $\cup H_{i}$ by taking as strata the
intersections of the hyperplanes, one stratum for each element of the
intersection lattice. \ This stratification is ``bonne'' in the sense of
\cite{Le1}, and $L$ is transverse to the strata, except $\{0\}$, since $L$
contains no other intersection of hyperplanes. \ The condition on $\Gamma$ is
required for L\^{e}'s result to ensure that $\left|  z_{\ell}\right|  _{F}$ is
Morse. \ See \cite{Randell} for further discussion.
\end{proof}

Such linear forms will be called generic. \ Note that such a generic form
yields an \emph{equivariant} Morse function (with respect to the covering
transformations) on $F$. \ (This group of covering transformations is
generated by $\zeta=\exp(2\pi i/n)$ acting by $\zeta(z_{1},z_{2},...,z_{\ell
})=(\zeta z_{1},\zeta z_{2},...,\zeta z_{\ell})$.) \ Therefore, as pointed out
by A. Dimca, one obtains a Morse function on the quotient $M^{\ast}$.

\begin{theorem}
For $L$ generic, the homotopy type of $F$ (resp. $M^{\ast}$) is obtained by
adding a number of $(\ell-1)$-cells to $F\cap\{L=0\}$ (resp. to $M^{\ast}%
\cap\{L=0\}$). \ The number of cells added to obtain $F$ is $n$ times the
number added to obtain $M^{\ast}$.
\end{theorem}

\begin{proof}
Clearly the number of critical points of the Morse function $\left|  L\right|
$ $=\left|  z_{\ell}\right|  $ on $F$ is $n$ times the number of critical
points of the induced Morse function on $M^{\ast}$. \ The Morse functions on
$F$ and $M^{\ast}$ are not proper, so we must take care. \ We follow here the
argument given by L\^{e} in \cite{Le1} for the Milnor fiber $F$. \ The same
argument works for $M^{\ast}$; one simply projects the sets chosen for $F$ via
the map $p:F\rightarrow M^{\ast}$.

The actual Milnor fiber is $F=Q^{-1}(1)$. \ We first note that for a
sufficiently large $r>>0$, $F$ is diffeomorphic to $F_{r}=F\cap\{z_{\ell}<r\}$
and homotopy equivalent to $\overline{F_{r}}=F\cap\{z_{\ell}\leq r\}$. \ Then
for $R$ sufficiently large with respect to $r$, we let $B_{R}$ be the closed
ball of radius $R$ in the first $(\ell-1)$-coordinates and note that $F$ has
the homotopy type of $F\cap(B_{R}\times D_{r})$, where $D_{r}$ is the set of
complex numbers $z_{\ell}<r.$ \ Further, from \cite[Proposition 2.1]{Le1}, we
note that $D_{r}$ and $B_{R}$ can be chosen so that for all $\eta\in
D_{r}-\{0\}$, $F\cap\{z_{\ell}=\eta\}$ is transverse to $\partial B_{R}%
\times\{\eta\}$ in $\{z_{\ell}=\eta\}$.

Then it follows that $\left|  z_{\ell}\right|  :F\cap(B_{R}\times
D_{r})\rightarrow D_{r}$ is a proper relative Morse function on a manifold
with boundary, and the map on the boundary is a submersion. \ Therefore
$F\cap(B_{R}\times D_{r})$ has the homotopy type of $F\cap\{z_{\ell}=0\}$ with
$(\ell-1)$-cells attached. \ As noted, the same argument works in $M^{\ast}$
by considering images under the map $p$.
\end{proof}

\ Here is a simple example.

\textbf{Example } \ Consider the $3$-arrangement given by $Q=x(x-y)(x+y-z)$
and consider the form $L=z$. \ Note that $\left\{  z=0\right\}  $ contains no
intersection of the hyperplanes. \ One computes the polar curve $\Gamma
=\{z-2y=0,3x-y=0\}$ and thus $\left|  L\right|  _{F}$ is a relative Morse
function on $F=Q^{-1}(1)$. \ The critical points are the solutions of the
system of equations
\[%
\begin{array}
[c]{c}%
Q=1\\
Q_{x}=0\\
Q_{y}=0
\end{array}
\]
Here as usual, $Q_{x}=\frac{\partial Q}{\partial x}$. \ \ One may calculate
directly that there are three critical points, each of index two. \ Thus $F$
is formed by adding three 2-cells to $F\cap\left\{  z=0\right\}  $. \ The set
$F\cap\left\{  z=0\right\}  $ is the Milnor fiber for a 2-arrangement, so we
repeat the process, this time taking $L=y$. \ This time one calculates that
there are six critical points of index one, so that $F\cap\left\{
z=0\right\}  $ is formed by adding six 1-cells to $F\cap\left\{
y=z=0\right\}  .$ \ Since this latter set is simply the solutions of $x^{3}%
=1$, or three points, we obtain a CW structure for the homotopy type of $F$
which has three $0$-cells, six $1$-cells, and three $2$-cells. \ Similarly, we
obtain a CW structure on the homotopy type of $M^{\ast}$ which has one
$0$-cell, two $1$-cells, and one $2$-cell (and of course $M^{\ast}$ has the
homotopy type of the two-torus). \ Finally, since $M\cong M^{\ast}%
\times\mathbb{C}^{\ast}$ we obtain a CW complex with the homotopy type of $M$
by using the product structure. \ The structures thus obtained on $M$ and
$M^{\ast}$ in this case (and always) have the number of cells in each
dimension equal to the rank of the homology group in that dimension, and in
this sense they are minimal structures. \ A more detailed explicit analysis
yields the attaching maps, and one can see that $F$ is a two-torus with a
rather inefficient cell structure.

Of course in this and what follows we could phrase results in terms of adding
handles, and thus build the diffeomorphism type of the Milnor fiber, not
simply the homotopy type. \ 

Next we ask how many critical points (handles, cells added) there are at each
stage. \ Let $b_{p}$ denote the usual Betti number.

\begin{theorem}
\label{Numbertheorem}The number of cells of dimension $p$ in the CW structure
for $M^{\ast}$ equals $b_{p}(M^{\ast}).$ \ The number of cells of dimension
$p$ in the CW structure for $F$ is $n\cdot b_{p}(M^{\ast})$. \ 
\end{theorem}

\begin{proof}
We use induction on $\ell$. \ For $\ell=1$ the result follows since $F$
consists of $n$ points, $M^{\ast}$ is a point, and $M\cong M^{\ast}%
\times\mathbb{C}^{\ast}$. \ For $\ell=2$ we note that $M$ is the complement of
a union of complex dimension one subspaces, so that $M$ has the homotopy type
of the complement of an $(n,n)$-torus link and $M^{\ast}$ is an n-times
punctured two-sphere. \ Thus $P(M)=(1+t)(1+(n-1)t)$. \ Observe that $F$ is
formed by attaching $c_{1}$ one-cells to $n$ points. \ To compute $c_{1}$ we
note that the Euler-Poincar\'{e} characteristic $e(M^{\ast})$ is $2-n$ so that
$e(F)=ne(M^{\ast})=n(2-n)=n-c_{1}.$ \ Thus $c_{1}=n(n-1)$, which is indeed $n$
times the first betti number of $M^{\ast}$. \ The case $\ell=2$ also follows
from the general induction step below.

For the general inductive step, we consider the $\ell-$arrangement
$\mathcal{A}$ with complement $M$. \ We choose a generic hyperplane $H$, and
consider the arrangement $\mathcal{A}^{H}=\{H\cap H_{i}\}$ as an arrangement
in the $(\ell-1)-$dimensional vector space $H$. \ Since $H$ is generic, no
intersection of hyperplanes from the original arrangement (a flat in the
language of arrangement theory) is contained in $H$ (except for \{0\}). \ Thus
the lattice of $\mathcal{A}^{H}$ is the same as that of the original
arrangement $\mathcal{A}$ in ranks (= codimension) $0$ through $\ell-2$.
\ Thus the betti numbers of $M$ are the same as those of the complement
$M^{H}$ of the union of the $H\cap H_{i}$ in $H$ for rank $0$ through
$\ \ell-2$ \cite{OrlikTeraoBook}. \ Since the Poincar\'{e} polynomials satisfy
$P(M^{\ast})=P(M)/(1+t)$, it follows that $b_{k}(M^{\ast})=b_{k}(M^{H\ast}) $
for $k\leq\ell-2$. \ We are attaching $(\ell-1)$-cells to $M^{H\ast} $ to form
$M^{\ast}$. \ These cells can only change homology groups in degree $(\ell-2)$
or $(\ell-1).$ \ But since $b_{\ell-2}(M^{\ast})=b_{\ell-2}(M^{H\ast})$ all
the $(\ell-1)$-cells attached to $M^{H\ast} $ to form $M^{\ast}$ are added
with homologically trivial attaching maps, so that the number of $(\ell
-1)$-cells added is $b_{\ell-1}(M^{\ast})$, as desired.
\end{proof}

Since $M$ has the homotopy type of $M^{\ast}\times S^{1}$, we form the product
CW structure for the homotopy type of $M$, taking the CW structure on $S^{1}$
consisting of a single $0$-cell and a single $1$-cell. \ Since the
Poincar\'{e} polynomial of $M$ is determined combinatorially by the lattice
and $P(M)=P(M^{\ast})(1+t)$ , one obtains a combinatorial determination of the
number of cells in the CW-structure for the homotopy types of $M$, $M^{\ast}$,
and $F$. \ \ \ We next state a number of corollaries.

\begin{corollary}
$M$ and $M^{\ast}$ have the homotopy type of CW complexes which have the
minimal number of cells in each dimension. \ That is, the number of $p$-cells
equals the $p$-th betti number.

\begin{definition}
(Papadima-Suciu) A space $X$ is \emph{minimal }provided that

i) $X$ is homotopy equivalent to a connected finite-type CW complex in which
the number of $p$-cells equals the $p$-th betti number of $X$.

ii) The homology groups of $X$ are torsion-free, and

iii) The cup-product map $\cup:$ $\wedge^{\ast}H^{1}(X)\rightarrow H^{\ast
}(X)$ is surjective.
\end{definition}
\end{corollary}

\begin{corollary}
$M$ and $M^{\ast}$ are minimal.
\end{corollary}

This result is also obtained by Dimca \cite{DimcaMilnor}.

\textbf{Example }\ It is not the case that complements of arbitrary
hypersurfaces are minimal. \ For example, the complement of the plane cusp
$\{x^{2}=y^{3}\}$ deformation retracts to the complement of the trefoil knot
in the three-sphere. \ As observed by \cite{PapaSucmin} any non-trivial knot
has a non-minimal complement. \ To see this, simply note that the first
homology group is infinite cyclic (so that $b_{1}=1$). \ If the complement had
a CW-complex structure with one cell of dimension one, then the fundamental
group would be cyclic. \ But no non-trivial knot complement has a cyclic
fundamental group.

From the above corollary one obtains from \cite[Theorem 1.4]{PapaSucmin} the
following two results. \ We let $\pi$ denote the fundamental group of $M$. \ 

\begin{corollary}
Suppose $Y=K(\pi,1)$ is a finite-type minimal CW-complex, with torsion free
cohomology, generated in degree 1$.$ \ Then if $b_{k}(Y)=b_{k}(M)$ for all
$k\leq p$,
\[
\pi_{k}(M)=0\text{, for }1<k<p\text{.}%
\]
\end{corollary}

Recall that a space $X$ is a $K(G,1)\,$\ if and only if its fundamental group
is isomorphic to $G$ and its universal cover is contractible. \ It is natural
to make the following

\begin{definition}
A group $G$ is \emph{minimal} if and only if there exists a minimal $K(G,1)$ space
\end{definition}

A well-known problem in the theory of arrangements is to determine under what
conditions the complement is a $K(\pi,1)$ space. \ A special case of the above
corollary is

\begin{corollary}
\label{iff}The complement $M$ of a complex hyperplane arrangement is a
$K(\pi,1)$ space if and only if $\pi$ is minimal and $b_{k}(K(\pi
,1))=b_{k}(M)$ for all $k$.
\end{corollary}

Thus whether or not $M$ is a $K(\pi,1)$ is determined by the group $\pi$ and
the lattice. \ The hypothesis on $\pi$ is necessary, as examples
($Q=xyz(y+z)(x-z)(2x+y)$) of arrangements are known in which the group
homology of the fundamental group is not finitely generated in degree three
\cite{Arvola}. \ It is known \cite[Theorem 13]{Randell} that it is always true
that $H^{i}(M)=H^{i}(\pi)$ for $i=1,2$.

\begin{corollary}
The number of cells in each dimension in this CW structure for $F,M,$ and
$M^{\ast}$ is determined by the combinatorics of the arrangement (i.e. the
intersection lattice.)
\end{corollary}

\ Artal Bartolo \cite{Artal-Bartolo} has found arrangements with lattices
yielding the same betti numbers for $M$, but different homology groups in the
Milnor fibers. \ Thus the attaching maps for the cells must differ.

\begin{corollary}
Letting $c_{p}$ denote the number of cells in the CW structure as found above
for $F$, and setting $C(F)=\sum c_{i}(F)t^{i},$ one has
\[
C(F)=nP(M^{\ast})=nP(M)/(1+t)
\]
\end{corollary}

In the example above, $P(M)=(1+t)^{3}$ and $n=3$ .

Finally, we examine more closely the location and number of critical points of
the Morse function $\left|  z\right|  $ on $F$. \ We discuss the situation for
$\ell=3$ ; the general situation is similar. \ Now the critical points are the
solutions of the system
\begin{equation}%
\begin{array}
[c]{c}%
Q=1\\
Q_{x}=0\\
Q_{y}=0
\end{array}
\label{systemaffine}%
\end{equation}
in $\mathbb{C}^{3}$ with variables $x,y,z$. \ We consider the related system
of equations
\begin{equation}%
\begin{array}
[c]{c}%
Q=w^{n}\\
Q_{x}=0\\
Q_{y}=0
\end{array}
\label{systemproj}%
\end{equation}
in the complex projective space $\mathbb{CP}^{3}$, with homogeneous
coordinates $[x:y:z:w]$. \ These equations have degree $n,n-1,n-1$
respectively, so that there are $n(n-1)^{2}$ solutions in projective space
(counting multiplicities). \ The critical points are just the solutions with
$w\neq0$.

Now consider a generic arrangement $\mathcal{B}$; that is, no three
hyperplanes intersect in a one-dimensional subspace of $\mathbb{C}^{3}$.
\ Then one may easily calculate that the CW structure on the associated Milnor
fiber has $n$ zero-cells, $n(n-1)$ one-cells, and $(n(n-1)(n-2))/2$ two-cells.
\ As a check, we calculate the number of solutions of (\ref{systemproj}) with
$w=0$. \ In this case, $Q=0,$ where $Q$ is the defining equation for the
arrangement $\mathcal{B}.$ \ Since $Q=\prod\alpha_{i}$ , this can only happen
when some $\alpha_{i}=0$. \ Now one easily computes $Q_{x}$ and $Q_{y},$%
\[
Q_{x}=\sum_{i=1}^{n}\frac{\partial\alpha_{i}}{\partial x}\frac{Q}{\alpha_{i}}%
\]
and
\[
Q_{y}=\sum_{i=1}^{n}\frac{\partial\alpha_{i}}{\partial y}\frac{Q}{\alpha_{i}}%
\]
Therefore when $\alpha_{i}=0$, $Q_{x}=0$ if and only if either $a_{j}=0$, some
$j\neq i$, or $\frac{\partial\alpha_{i}}{\partial x}=0$. \ Similarly, if
$\alpha_{i}=0$ , then $Q_{y}=0$ if and only if either $a_{j}=0$, some $j\neq
i$, or $\frac{\partial\alpha_{i}}{\partial y}=0$. \ But if both partials
$\frac{\partial\alpha_{i}\text{ }}{\partial x}$ and $\frac{\partial\alpha_{i}%
}{\partial y}$ are simultaneously zero, then $\alpha_{i}=z$, which does not
hold since we chose the Morse function as $\left|  z\right|  ,$ where the
hyperplane $z=0$ contains no flats of the arrangement. \ Therefore, solutions
of (\ref{systemproj}) with $w=0$ occur exactly at points (of the
$\mathbb{CP}^{2}$ where $\ w=0$) where two $a_{i}$ are zero. \ Thus, for
generic $\mathcal{B}$, there are $\binom{n}{2}$ such points, each with
multiplicity $n$ (the multiplicity of $Q=0$), or $n\binom{n}{2}$ solutions
counted with multiplicities. \ Thus there are
\[
n(n-1)^{2}-n\binom{n}{2}=\frac{n(n-1)(n-2)}{2}%
\]
solutions of (\ref{systemaffine}), and thus the expected number of critical
points. \ For non-generic $3-$arrangements, there are fewer $2$-cells added,
hence fewer critical points, hence \emph{more }\ solutions of
(\ref{systemproj}) with $w=0$. \ In the example below we will observe
solutions tending to $w=0$ in degenerating families.

The above argument shows in general (any arrangement, any dimension) that the
solutions of (\ref{systemproj}) with $w=0$ occur with two $\alpha_{i}$ equal
to zero. \ Therefore, we have

\begin{proposition}
The number of critical points, and hence the number of $(\ell-1)-$cells added
to form $F$, is given by
\[
n(n-1)^{\ell-1}-\left|  \{\text{solutions, with multiplicities, of
(}\ref{systemproj}\text{) on }\{w=0\}\}\right|
\]
\end{proposition}

\textbf{Example }\ Let $Q_{t}=x(x+y)(x-y+tz).$ \ With $t\neq0$, the three
hyperplanes intersect only at the origin, while for $t=0$, they intersect
along the $z$-axis. \ That is, at the point $[0:0:1]$ \ in homogeneous
coordinates for $M^{\ast}\subset CP^{2}$. \ Our earlier example was in fact
the case $t=1;$ all non-zero values for $t$ exhibit the same behavior. \ But
if $t=0$, the Morse function $\left|  z\right|  $ has no critical points (and
the polar curve $\Gamma$ is empty). \ To see this, note that at a critical
point of $\left|  z\right|  $ one must have
\[%
\begin{array}
[c]{c}%
Q=x(x+y)(x-y)=x^{3}-xy^{2}=1\\
Q_{x}=(x+y)(x-y)+x(x-y)+x(x+y)=3x^{2}-y^{2}=0\\
Q_{y}=x(x-y)-x(x+y)=-2xy=0
\end{array}
\]
The third equation shows that $x=0$ or $y=0$. \ But $x=0$ is clearly
impossible, from the first equation, while $y=0$ forces $x=0$ in the second
equation, again impossible by the first equation. \ So in fact $F_{0}$ has the
homotopy type of a 1-complex, obtained by removing (or not adding) the three
two-cells required to form $F_{1}$.


\begin{thebibliography}{9}
\bibitem{Milnor1}Milnor, J., Singular points of complex hypersurfaces, Annals
of Math. Studies \textbf{61}, Princeton University Press, 1968.

\bibitem{Le1}L\^{e}, D. T.,Calcul du nombre de cycles \'{e}vanouissants d'une
hypersurface complexe, Ann. Inst. Fourier, Grenoble \textbf{23}, 4 (1973), 261-270.

\bibitem{PapaSucmin}Papadima, S. and Suciu, A., Higher homotopy groups of
complements of complex hyperplane arrangements, preprint, arXiv:math.AT/0002251

\bibitem{OrlikTeraoBook}Orlik, P., and Terao, H., Arrangements of Hyperplanes,
Grundlehren der mathematischen Wissenschaften \textbf{300}, Springer Verlag, 1992.

\bibitem{OTMorse}Orlik, P. and Terao, H. Arrangements and Milnor fibers, Math.
Ann. \textbf{301}, (1995), 211-235.

\bibitem{DimcaMilnor}Dimca, A., Hypersurface complements, Milnor fibers and
minimality of arrangements, preprint, arXiv:math.AG/0011222.

\bibitem{Artal-Bartolo}Artal Bartolo, E. , Combinatorics and topology of line
arrangements in the complex projective plane, Proc. Amer. Math. Soc.,
\textbf{121}, (1994), 385-390.

\bibitem{Arvola}Arvola, B., Arrangements and cohomology of groups, preprint, 1992.

\bibitem{Randell}Randell, R., Homotopy and group cohomology of arrangements,
Topology and its Applications \textbf{78, }(1997), 201-213.
\end{thebibliography}
\end{document}